Article CJB/2010/7

# Generalization of the Wallace-Simson line

**Christopher J Bradley**

**Abstract:** A Wallace-Simson line of a point *J* has the property that it passes through the midpoint *Q* of *HJ*, where *H* is the orthocentre. By allowing *Q* to move on the perpendicular bisector of *HJ* we obtain generalized Wallace-Simson lines. These turn out to be known in another context.

## 1. Introduction

If *ABC* is a triangle and *J* is a point on its circumcircle *Σ*, then the Wallace-Simson line property is that the feet of the perpendiculars *L*, *M*, *N* from *J* on to the sides of *ABC* are collinear and *LMN* is called the Wallace-Simson line. Furthermore, if *H* is the orthocentre of *ABC* then the midpoint *Q* of *JH* lies on *LMN*. If instead of taking the feet of the perpendiculars on to the sides, we take the reflections *L'*, *M'*, *N'* of *J* in the sides, then *L'M'N'* is also a straight line which we refer to as the *Double Wallace-Simson line* of *J*. As this line is a degenerate Hagge circle it is known that *L'M'N'* passes through the orthocentre *H*.

In this article we provide a generalization of this construction, so that the Wallace-Simson line is just a particular case, being one of an infinite number of lines that we refer to as generalized Wallace-Simson lines. Later we identify these lines as ones that are already known, providing an explanation of why these lines are the same as those we describe.

Fig. 1 illustrates the general construction. The point *Q* is now varied from its fixed position at the midpoint of *JH* from which the well known Wallace-Simson line results and becomes a variable point on the perpendicular bisector of *JH*. The circle *S* is drawn, centre *Q*, passing through *J* and *H*. The points where the altitudes of *ABC* meet *S* are denoted by *X*, *Y*, *Z*, where *AX*, *BY*, *CZ* are perpendicular to *BC*, *CA*, *AB* respectively. Circles *AJX*, *BJY*, *CJZ* are then drawn and *L* is the point of intersection of *BJY* and *CJZ* with *M*, *N* similarly defined. As we prove in due course, the following results hold:

1. *L*, *M*, *N* lie on *BC*, *CA*, *AB* respectively;
2. *L*, *M*, *N* are collinear.
3. The line *LMN* contains the point *Q*.

The line *LMN* is, of course, what we mean by a generalized Wallace-Simson line and as *Q* varies we get an infinite number of such lines all defined with respect to the one fixed point *J* on the circumcircle of *ABC*.

Various subsidiary results emerge in the course of establishing these propositions, and we now briefly indicate the plan we adopt to prove them.

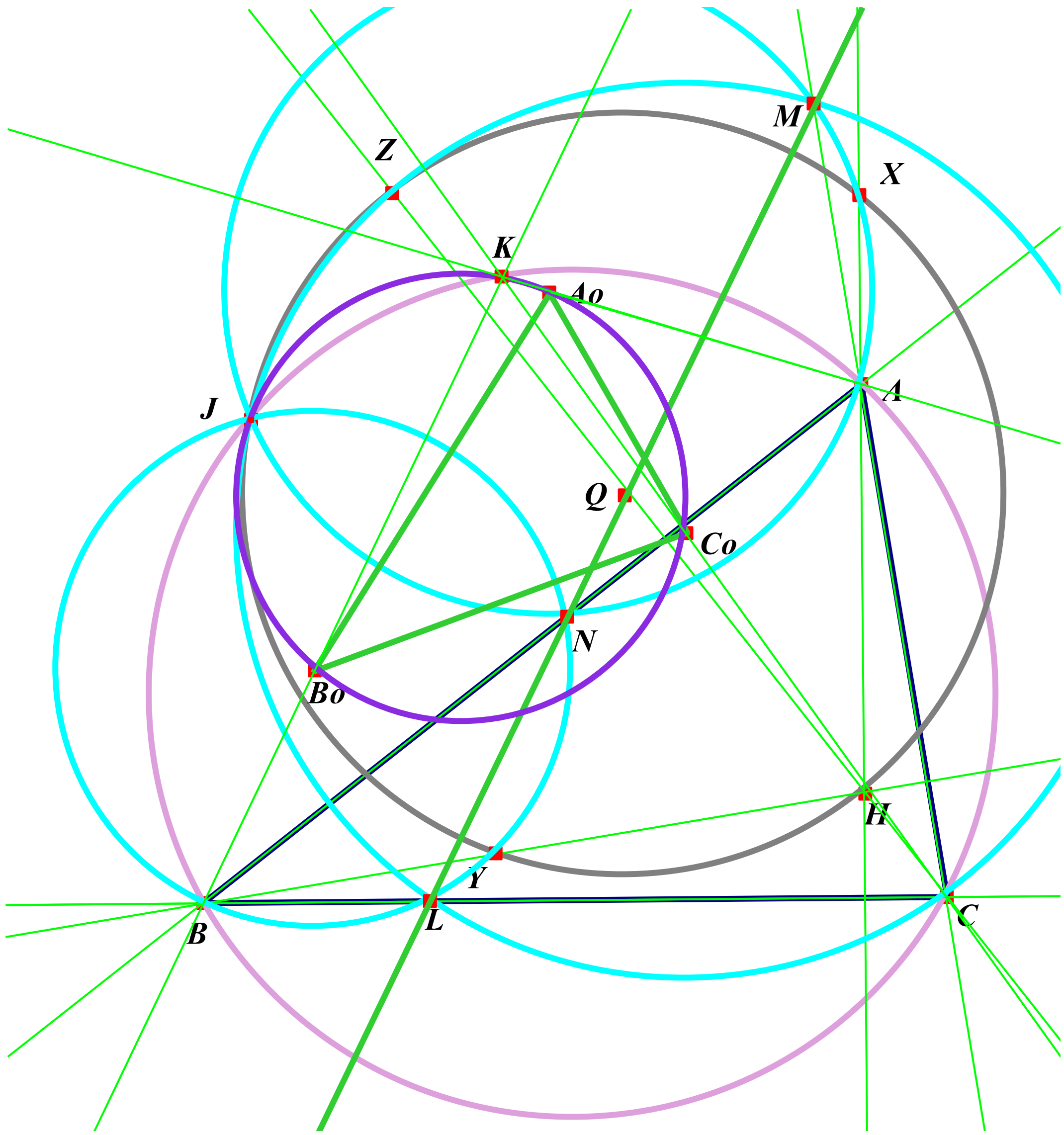

**Fig. 1**
**The generalized Wallace-Simson line LMN**

The method we use is to take *J* to be the centre of the direct similarity that maps *H* to $Q$. We then apply this direct similarity to $ABC$ to get a triangle $A_oB_oC_o$ and it transpires that $\Sigma_o$, the circumcircle of $A_oB_oC_o$, shares the point $J$ with $\Sigma$ and a second point $K$ such that $AA_o$, $BB_o$, $CC_o$ concur at $K$. This is typical of directly similar triangles in perspective, see Wood [1]. When $Q$ is the midpoint of $JH$ triangle $A_oB_oC_o$ has its sides parallel to $ABC$ and is half the size. In all cases it turns out that the line $LMN$ we want is the Double Wallace- Simson line of $J$ with respect to triangle $A_oB_oC_o$ and, being a Double Wallace-Simson, line it automatically passes through the orthocentre of $A_oB_oC_o$, which by construction is $Q$. Synthetic considerations may be used to show that $L$, $M$, $N$ lie on $BC$, $CA$, $AB$ respectively and indeed that $LMN$ is a straight line. The bonus of the analysis is that it highlights the role played by the direct similarity.

## 2. Defining the direct similarity

We use Cartesian co-ordinates and take $J$ to be the origin (0, 0) and the centre $O$ of triangle $ABC$ to the point (1, 0). Throughout the calculation the computer algebra package *DERIVE* was used. If, for the moment, we suppose $H$ has co-ordinates $(h, k)$, then it is straightforward to show that a variable point $Q$ on the perpendicular bisector of $JH$ has co-ordinates $(\tfrac{1}{2}h - kt, \tfrac{1}{2}k + ht)$, where $t$ is a variable parameter.

The direct similarity through $J$ that maps $H$ to $Q$ is therefore given by the matrix

$$\begin{pmatrix} 1/2 & -t \\ t & 1/2 \end{pmatrix}. \tag{2.1}$$

If considered as a dilation of magnitude $d$ with centre of enlargement $J$ and a rotation by an angle $\theta$ about $J$ we see that $d$ and $\theta$ are related by the equation $d = 1/(2\cos\theta)$ and also that $t = \tfrac{1}{2}\tan\theta$. (In all that follows we suppose without loss of generality that $\cos\theta$ is positive and that $-\tfrac{1}{2}\pi < \theta < \tfrac{1}{2}\pi$.)

The equation of the circumcircle $\Sigma$ is

$$x^2 + y^2 - 2x = 0. \tag{2.2}$$

We suppose that the point $A$ has co-ordinates $\left(\frac{2}{1+a^2}, \frac{2a}{1+a^2}\right)$, where $a$ is a parameter, and that $B$, $C$ have similar co-ordinates with $b$, $c$ respectively replacing $a$. Applying the direct similarity to $A$ we obtain a point $A_o$ with co-ordinates $\left(\frac{1-2at}{1+a^2}, \frac{a+2t}{1+a^2}\right)$. The points $B_o$, $C_o$, the images of $B$, $C$ respectively under the direct similarity have similar co-ordinates with $b$, $c$
respectively replacing $a$.

We now prove that triangle $A_oB_oC_o$ is in perspective with triangle $ABC$. The equation of $AA_o$ is

$$(a - 2t)x - (1 + 2at)y + 4t = 0. \qquad (2.3)$$

This meets the circumcircle $\Sigma$ with Equation (2.2) at $A$ and at the point $K$ with co-ordinates $\left(\frac{8t^2}{1+4t^2}, \frac{4t}{1+4t^2}\right)$. Since the co-ordinates of $K$ are independent of $a$, $b$, $c$ it follows that $AA_o$, $BB_o$, $CC_o$ are concurrent at $K$ and the triangles are not only similar, but in perspective. From the theory of Wood [1] on directly similar triangles in perspective it follows that the circle $\Sigma_o$, the circumcircle of triangle $A_oB_oC_o$ intersects $\Sigma$ at the points $J$ and $K$. See Fig. 1.

We now consider the circle, centre $A_o$ passing through $J$. Its equation is soon found to be

$$(1 + a^2)x^2 + (1 + a^2)y^2 - 2(1 - 2at)x - 2(a + 2t)y = 0. \qquad (2.4)$$

It may now be checked that this circle passes through $A$. So this circle is now identified as circle $JAX$, though we have yet to find the co-ordinates of $X$ and to show it lies on the Hagge circle, centre $Q$ passing through $H$ and $J$. If the parameter $a$ in Equation (2.4) is replaced by $b$, $c$ we have the equations of the circles that may be identified as circles $JBY$, $JCZ$ respectively. Note that in view of the direct similarity triangles $JA_oA$, $JB_oB$, $JC_oC$ are similar.

The co-ordinates of $H$ may be obtained from those of $O$, $A$, $B$, $C$ and are

$$\left(\frac{2(2+a^2+b^2+c^2-2a^2b^2c^2)}{(1+a^2)(1+b^2)(1+c^2)}, \frac{2(a+b+c+ab^2c^2+bc^2a^2+ca^2b^2+ab^2+ac^2+bc^2+ba^2+ca^2+cb^2)}{(1+a^2)(1+b^2)(1+c^2)}\right). \qquad (2.5)$$

From this we may obtain the equation of the altitude $AH$, which is

$$(1 + a^2)(b + c)x - (1 + a^2)(1 - bc)y + 2(a + b + c - abc) = 0. \qquad (2.6)$$

This meets the circle centre $A_o$ through $J$ at the point $A$ and the point $X$ with co-ordinates

$$\left(\frac{2(b+c+2t-2bct)(abc-a+b+c)}{(1+a^2)(1+b^2)(1+c^2)}, \frac{2(abc-a+b+c)(bc+2t(b+c)-1)}{(1+a^2)(1+b^2)(1+c^2)}\right). \qquad (2.7)$$

The co-ordinates of $Y$ and $Z$ may be found from those of $X$ by cyclic change of the symbols $a$, $b$, $c$. The equation of the circle through $X$, $Y$, $Z$ may now be calculated and is

(2.8)

$$(1 + a^2)(1 + b^2)(1 + c^2)(x^2 + y^2) + 2(a^2b^2c^2 + 2abct(bc + ca + ab) + 2t(a^2b + a^2c + b^2c + b^2a + c^2a + c^2b) + 2t(a + b + c) - a^2 - b^2 - c^2 - 2)x + 2(2a^2b^2c^2t - abc(bc + ca + ab) - 2t(a^2 + b^2 + c^2) - (a^2b + a^2c + b^2c + b^2a + c^2a + c^2b) - (a + b + c + 4t))y = 0.$$

This circle clearly passes through $J$ and it may be verified by substitution that it also passes through $H$. Equation (2.8) therefore represents the Hagge circle, centre $Q$ and passing through $J$ and $H$. Fig. 1 has now been completely analysed and from here the results stated in the introduction follow by synthetic arguments and other known results.

## 3. Results

The angle arguments given below are appropriate to the figure at the start of the article and as such are diagram dependent. However, they are easily modified for other possible figures, or may be written in a diagram independent way by working mod $\pi$.

**Theorem 4.1**
Angle $JLB$ = Angle $JYB$ = Angle $JYH$ = $180^{o}$ – Angle $JZH$ = $180^{o}$ – Angle $JZC$ = Angle $JLC$.
This shows that $L$ lies on $BC$. Similarly, $M$ and $N$ lie on $CA$ and $AB$ respectively.

**Theorem 4.2**
Now angle $JNL$ = Angle $JBL$ = = $180^{o}$ – Angle $JBC$ = Angle $JAC$ = Angle $JAM$ = $180^{o}$ – Angle $JNM$ and so $LMN$ is a straight line.

**Theorem 4.3**
Circles $JBY$ and $JCZ$ have centres $B_o$ and $C_o$ respectively. The common chord $JL$ is perpendicular to the line of centres $B_oC_o$ and bisected by it. $L$ is therefore the reflection of $J$ in the side $B_oC_o$ of triangle $A_oB_oC_o$. The same holds for $M$ and $N$ with respect to the sides $C_oA_o$ and $A_oB_o$ respectively. $LMN$ is therefore the Double Wallace-Simson line of $J$ with respect to triangle $A_oB_oC_o$. But, by construction, $Q$ is the orthocentre of this triangle and therefore lies on $LMN$.

Since $LMN$ is drawn making right angles with the sides of a triangle which is similar and at a fixed angle to the sides of $ABC$, it is clear that these lines make equal angles with the sides of $ABC$. It is interesting to observe that these generalized Wallace-Simson lines (sometimes called *oblique* Wallace-Simson lines) exist only because Double Wallace-Simson lines exist for a smaller similar triangle. The original proof seems to have been by Carnot, the best reference I could find being in Exercises by Frére Gabriel-Marie [2].

Flat 4,
Terrill Court,
12-14 Apsley Road,
BRISTOL BS8 2SP